\documentclass[12pt, a4paper]{article}
\usepackage[latin1]{inputenc}
\usepackage[english]{babel}
\usepackage{amsmath}
\usepackage[english]{babel}
\usepackage{indentfirst}
\usepackage{amstext}
\usepackage{enumerate}
\usepackage{amsfonts}
\usepackage{textcomp}
\usepackage{amssymb}
\usepackage[dvips]{graphicx}
\usepackage{setspace}
\usepackage[all]{xy}

\newcommand {\demo}{\hskip -0.6cm{\bf Proof.  }}
\newcommand {\fim}{\hfill{$\square$}\vskip 1pc}
\newcommand {\nl}{\newline}

\newcommand {\R}{\mathbb{R}}

\newcommand {\C}{\mathbb{C}}

\newcommand {\GG}{\mathcal{G}}
\newcommand {\FF}{\mathcal{F}}
\newcommand {\T}{\mathbb{T}}

\newcommand {\G}{\mathcal{G}}
\newcommand {\F}{\mathbb{F}}

\newcommand{\Aut}{\operatorname{Aut}}
\newcommand{\lsp}{\overline{\operatorname{span}}}

\newtheorem{teorema}{Theorem}[section]
\newtheorem{lema}[teorema]{Lemma}
\newtheorem{corolario}[teorema]{Corollary}
\newtheorem{definicao}[teorema]{Definition}
\newtheorem{proposicao}[teorema]{Proposition}

\newtheorem{remark}[teorema]{Remark}
\DeclareMathOperator{\aut}{Aut}

\begin{document}
\onehalfspace

\title{KMS and ground states on ultragraph C*-algebras}
\maketitle

\begin{center}
\author{Gilles Gon\c{c}alves de Castro and Daniel Gon\c{c}alves\footnote{Partially supported by Conselho Nacional de Desenvolvimento Científico e Tecnológico (CNPq) - Brazil}}
\end{center}

\vspace{0.5pc}

\begin{abstract}
We describe KMS and ground states arising from a generalized gauge action on ultragraph C*-algebras. We focus on ultragraphs that satisfy Condition~(RFUM), so that we can use the partial crossed product description of ultragraph C*-algebras recently described by the second author and Danilo Royer. In particular, for ultragraphs with no sinks, we generalize a recent result by Toke Carlsen and Nadia Larsen: Given a time evolution on the C*-algebra of an ultragraph, induced by a function on the edge set, we characterize the KMS states in five different ways and ground states in four different ways. In both cases we include a characterization given by maps on the set of generalized vertices of the ultragraph. We apply this last result to show the existence of KMS and ground states for the ultragraph C*-algebra that is neither an Exel-Laca nor a graph C*-algebra.
\end{abstract}

\vspace{0.5pc}



\section{Introduction}

KMS (Kubo-Martin-Schwinger) states on C*-algebras have been the subject of intense research both in Mathematics and Physics. The mathematical study of KMS states stemmed from the realization that in the C*-algebraic formulation of quantum statistical mechanics, KMS states are equilibrium states associated to the one-parameter group of automorphisms of time evolution (see \cite{MR1441540}). Over the years a large literature was developed, as researches studied KMS states on C*-algebras associated to expansive maps (see \cite{MR2215776}), graphs (see \cite{MR3088995, MR3306919}), relative graphs (see \cite{MR3539347}), higher rank graphs (see \cite{MR3121730}), local homemorphisms (see \cite{MR3254419}), and Fell bundles over groupoids (see \cite{1708.00629}), to name a few. It is interesting to note that, as mentioned in \cite{1708.00629}, KMS states make sense for any C*-dynamical system, and there is significant evidence that the KMS data is a useful invariant of a dynamical system. For example, the results in \cite{MR759450} show that the KMS data for a Cuntz-Krieger algebra encodes the topological entropy of the associated shift space.

Ultragraphs are generalizations of graphs, where the range map take values on the power set of the vertices. Ultragraph C*-algebras were introduced by Tomforde in \cite{MR2050134} and they were key in the study of the relations between Exel-Laca and graph algebras (see \cite{MR2629692}), provided new examples of algebras associated to combinatorial data (see \cite{1701.00323}, \cite{MR2001938}), and are connected to the theory of the Perron-Frobenious operator (see \cite{MR3554458}). Recently ultragraphs and their C*-algebras have played an important role in the study of shift spaces associated to infinite alphabets, see \cite{MR3600124, GRultrapartial}. In particular, based on ultragraph C*-algebra theory, in \cite{GRultrapartial} a new approach to shifts of finite type over infinite alphabets was proposed. 

Given the above mentioned relation between KMS data and topological entropy of shifts, it is relevant to study KMS states associated to ultragraph C*-algebras. This will be the focus of this article. In our main theorem (see Theorem~\ref{teoremaprincipal}), given an ultragraph $\G$ with no sinks that satisfy Condition~(RFUM), we identify the KMS states on $C^*(\G)$ with a subset of the states on a Abelian algebra ($C(X)$), a subset of the regular Borel probability measures on $X$, a subset of the functions on the generalized vertices of $\G$, and a subset of the states on the core subalgebra of $C^*(\G)$. Taking advantage of the results proved to characterize the KMS states we also describe the ground states on $C^*(\G)$.

Our work is organized in the following way: In Section~\ref{sec1} we review the construction of ultragraph C*-algebras via partial crossed products given in \cite{GRultrapartial} and prove a few auxiliary results. In particular we give a simpler description of the basis of the topological (shift) space associated to an ultragraph in \cite{GRultrapartial}. In Section~3 we study arbitrary ultragraphs $\G$ and show that there exists a correspondence between the KMS of the C*-algebra $C^*(\G)$ and a subset of the states on the core subalgebra of $C^*(\G)$. In order to use the results in \cite{MR1953065}, regarding KMS states of partial dynamical systems, in Section~4 we restrict ourselves to ultragraphs with no sinks that satisfy Condition~(RFUM) and prove Theorem~\ref{teoremaprincipal}, which generalizes (for ultragraphs with no sinks) Theorem~4.1 in \cite{MR3539347} (we remark that our methods to prove this result are in great part different from the methods in \cite{MR3539347}). In Section~\ref{ground} we describe ground states in terms of subsets of the states on an Abelian algebra, a subset of the regular Borel probability measure on the spectrum of the just mentioned Abelian algebra, and as a subset of the functions on the generalized vertices of $\G$. We end the paper in Section~\ref{secex}, where we show the existence of KMS for the ultragraph C*-algebra that is neither an Exel Laca nor a graph algebra and completely characterize the ground states on the aforementioned algebra.

\section{Ultragraph C*-algebras viewed as partial crossed products}\label{sec1}

\subsection{Ultragraph C*-algebras}

In this short subsection we recall the definition of ultragraph C*-algebras and a few key results, as done in \cite{MR2050134}.

\begin{definicao}\label{def of ultragraph}
An \emph{ultragraph} is a quadruple $\mathcal{G}=(G^0, \mathcal{G}^1, r,s)$ consisting of two countable sets $G^0, \mathcal{G}^1$, a map $s:\mathcal{G}^1 \to G^0$, and a map $r:\mathcal{G}^1 \to P(G^0)\setminus \{\emptyset\}$, where $P(G^0)$ stands for the power set of $G^0$.
\end{definicao}

Before we define the C*-algebra associated to an ultragraph we need the following notion.

\begin{definicao}\label{def of mathcal{G}^0}
Let $\mathcal{G}$ be an ultragraph. Define $\mathcal{G}^0$ to be the smallest subset of $P(G^0)$ that contains $\{v\}$ for all $v\in G^0$, contains $r(e)$ for all $e\in \mathcal{G}^1$, and is closed under finite unions and nonempty finite intersections.
\end{definicao}

The following description of $\GG^0$ is useful.

\begin{lema}\label{description}\cite[Lemma~2.12]{MR2050134}
If $\mathcal{G}$ is an ultragraph, then \begin{align*} \mathcal{G}^0 = \{
\bigcap_{e
\in X_1} r(e)
\cup \ldots 
\cup \bigcap_{e \in X_n} r(e) \cup F : & \ \text{$X_1,
\ldots, X_n$ are finite subsets of $\mathcal{G}^1$} \\ & \text{ and $F$
is a finite subset of $G^0$} \}.
\end{align*}  
Furthermore, $F$ may be chosen to
be disjoint from $\bigcap_{e \in X_1} r(e) \cup \ldots 
\cup \bigcap_{e \in X_n} r(e)$.
\end{lema}

\begin{definicao}\label{def of C^*(mathcal{G})}
Let $\mathcal{G}$ be an ultragraph. The \emph{ultragraph algebra} $C^*(\mathcal{G})$ is the universal $C^*$-algebra generated by a family of partial isometries with orthogonal ranges $\{s_e:e\in \mathcal{G}^1\}$ and a family of projections $\{p_A:A\in \mathcal{G}^0\}$ satisfying
\begin{enumerate}
\item\label{p_Ap_B=p_{A cap B}}  $p_\emptyset=0,  p_Ap_B=p_{A\cap B},  p_{A\cup B}=p_A+p_B-p_{A\cap B}$, for all $A,B\in \mathcal{G}^0$;
\item\label{s_e^*s_e=p_{r(e)}}$s_e^*s_e=p_{r(e)}$, for all $e\in \mathcal{G}^1$;
\item $s_es_e^*\leq p_{s(e)}$ for all $e\in \mathcal{G}^1$; and
\item\label{CK-condition} $p_v=\sum\limits_{s(e)=v}s_es_e^*$ whenever $0<\vert s^{-1}(v)\vert< \infty$.
\end{enumerate}
\end{definicao}

\subsection{Notation}
Next we set up notation that will be used throughout the paper. This agrees with notation introduced in \cite{MR2457327} and \cite{MR3600124}.

Let $\mathcal{G}$ be an ultragraph. A \textit{finite path} in $\mathcal{G}$ is either an element of $\mathcal{G}%
^{0}$ or a sequence of edges $e_{1}\ldots e_{k}$ in $\mathcal{G}^{1}$ where
$s\left(  e_{i+1}\right)  \in r\left(  e_{i}\right)  $ for $1\leq i\leq k$. If
we write $\alpha=e_{1}\ldots e_{k}$, the length $\left|  \alpha\right|  $ of
$\alpha$ is just $k$. The length $|A|$ of a path $A\in\mathcal{G}^{0}$ is
zero. We define $r\left(  \alpha\right)  =r\left(  e_{k}\right)  $ and
$s\left(  \alpha\right)  =s\left(  e_{1}\right)  $. For $A\in\mathcal{G}^{0}$,
we set $r\left(  A\right)  =A=s\left(  A\right)  $. The set of
finite paths in $\mathcal{G}$ is denoted by $\mathcal{G}^{\ast}$. An \textit{infinite path} in $\mathcal{G}$ is an infinite sequence of edges $\gamma=e_{1}e_{2}\ldots$ in $\prod \mathcal{G}^{1}$, where
$s\left(  e_{i+1}\right)  \in r\left(  e_{i}\right)  $ for all $i$. The set of
infinite paths  in $\mathcal{G}$ is denoted by $\mathfrak
{p}^{\infty}$. The length $\left|  \gamma\right|  $ of $\gamma\in\mathfrak
{p}^{\infty}$ is defined to be $\infty$. A vertex $v$ in $\mathcal{G}$ is
called a \emph{sink} if $\left|  s^{-1}\left(  v\right)  \right|  =0$ and is
called an \emph{infinite emitter} if $\left|  s^{-1}\left(  v\right)  \right|
=\infty$. 

For $n\geq1,$ we define
$\mathfrak{p}^{n}:=\{\left(  \alpha,A\right)  :\alpha\in\mathcal{G}^{\ast
},\left\vert \alpha\right\vert =n,$ $A\in\mathcal{G}^{0},A\subseteq r\left(
\alpha\right)  \}$. We specify that $\left(  \alpha,A\right)  =(\beta,B)$ if
and only if $\alpha=\beta$ and $A=B$. We set $\mathfrak{p}^{0}:=\mathcal{G}%
^{0}$ and we let $\mathfrak{p}:=\coprod\limits_{n\geq0}\mathfrak{p}^{n}$. We embed the set of finite paths $\GG^*$ in $\mathfrak{p}$ by sending $\alpha$ to $(\alpha, r(\alpha))$. We
define the length of a pair $\left(  \alpha,A\right)  $, $\left\vert \left(
\alpha,A\right)  \right\vert $, to be the length of $\alpha$, $\left\vert
\alpha\right\vert $. We call $\mathfrak{p}$ the \emph{ultrapath space}
associated with $\mathcal{G}$ and the elements of $\mathfrak{p}$ are called
\emph{ultrapaths}. Each $A\in\mathcal{G}^{0}$ is regarded as an ultrapath of length zero and can be identified with the pair $(A,A)$. We may extend the range map $r$ and the source map $s$ to
$\mathfrak{p}$ by the formulas, $r\left(  \left(  \alpha,A\right)  \right)
=A$, $s\left(  \left(  \alpha,A\right)  \right)  =s\left(  \alpha\right)
$ and $r\left(  A\right)  =s\left(  A\right)  =A$.

We concatenate elements in $\mathfrak{p}$ in the following way: If $x=(\alpha,A)$ and $y=(\beta,B)$, with $|x|\geq 1, |y|\geq 1$, then $x\cdot y$ is defined if and only if
$s(\beta)\in A$, and in this case, $x\cdot y:=(\alpha\beta,B)$. Also we
specify that:
\begin{equation}
x\cdot y=\left\{
\begin{array}
[c]{ll}%
x\cap y & \text{if }x,y\in\mathcal{G}^{0}\text{ and if }x\cap y\neq\emptyset\\
y & \text{if }x\in\mathcal{G}^{0}\text{, }\left|  y\right|  \geq1\text{, and
if }x\cap s\left(  y\right)  \neq\emptyset\\
x_{y} & \text{if }y\in\mathcal{G}^{0}\text{, }\left|  x\right|  \geq1\text{,
and if }r\left(  x\right)  \cap y\neq\emptyset
\end{array}
\right.  \label{specify}%
\end{equation}
where, if $x=\left(  \alpha,A\right)  $, $\left|  \alpha\right|  \geq1$ and if
$y\in\mathcal{G}^{0}$, the expression $x_{y}$ is defined to be $\left(
\alpha,A\cap y\right)  $. Given $x,y\in\mathfrak{p}$, we say that $x$ has $y$ as an initial segment if
$x=y\cdot x^{\prime}$, for some $x^{\prime}\in\mathfrak{p}$, with $s\left(
x^{\prime}\right)  \cap r\left(  y\right)  \neq\emptyset$. 

We extend the source map $s$ to $\mathfrak
{p}^{\infty}$, by defining $s(\gamma)=s\left(  e_{1}\right)  $, where
$\gamma=e_{1}e_{2}\ldots$. We may concatenate pairs in $\mathfrak{p}$, with
infinite paths in $\mathfrak{p}^{\infty}$ as follows. If $y=\left(
\alpha,A\right)  \in\mathfrak{p}$, and if $\gamma=e_{1}e_{2}\ldots\in
\mathfrak{p}^{\infty}$ are such that $s\left(  \gamma\right)  \in r\left(
y\right)  =A$, then the expression $y\cdot\gamma$ is defined to be
$\alpha\gamma=\alpha e_{1}e_{2}...\in\mathfrak{p}^{\infty}$. If $y=$
$A\in\mathcal{G}^{0}$, we define $y\cdot\gamma=A\cdot\gamma=\gamma$ whenever
$s\left(  \gamma\right)  \in A$. Of course $y\cdot\gamma$ is not defined if
$s\left(  \gamma\right)  \notin r\left(  y\right)  =A$. 

\begin{remark} When no confusion arises we will omit the dot in the notation of concatenation defined above, so that $x\cdot y$ will be denoted by $xy$.
\end{remark}

\begin{definicao}
\label{infinte emitter} For each subset $A$ of $G^{0}$, let
$\varepsilon\left(  A\right)  $ be the set $\{ e\in\mathcal{G}^{1}:s\left(
e\right)  \in A\}$. We shall say that a set $A$ in $\mathcal{G}^{0}$ is an
\emph{infinite emitter} whenever $\varepsilon\left(  A\right)  $ is infinite.
\end{definicao}



\subsection{The shift space associated to an ultragraph}

Our goal in this subsection is two-fold. We recall the shift space $X$ associated to an ultragraph without sinks (as in \cite{GRultrapartial}) and prove some new results about its topology. 

The key concept in the definition of the shift space $X$ associated to $\GG$ is that of minimal infinite emitters. We recall this below.

\begin{definicao}\label{minimal} Let $\GG$ be an ultragraph and $A\in \GG^0$. We say that $A$ is a minimal infinite emitter if it is an infinite emitter that contains no proper subsets (in $\GG^0$) that are infinite emitters. Equivalently, $A$ is a minimal infinite emitter if it is an infinite emitter and has the property that, if $B\in \GG^0$ is an infinite emitter, and $B\subseteq A$, then $B=A$. For a finite path $\alpha$ in $\GG$, we say that $A$ is a minimal infinite emitter in $r(\alpha)$ if $A$ is a minimal infinite emitter and $A\subseteq r(\alpha)$. We denote the set of all minimal infinite emitters in $r(\alpha)$ by $M_\alpha$.
\end{definicao}

For later use we recall the following Lemma.

\begin{lema}\label{miniftyem} Let $x=(\alpha, A) \in \mathfrak{p}$ and suppose that $A$ is a minimal infinite emitter. If the cardinality of $A$ is finite then it is equal to one and, if the cardinality of $A$ is infinite, then $A = \bigcap\limits_{e\in Y} r(e)$ for some finite set $Y\subseteq \GG^1$.
\end{lema}



Associated to an ultragraph with no sinks, we have the topological space $X= \mathfrak{p}^{\infty} \cup X_{fin}$, where 
$$X_{fin} = \{(\alpha,A)\in \mathfrak{p}: |\alpha|\geq 1 \text{ and } A\in M_\alpha \}\cup
 \{(A,A)\in \GG^0: A \text{ is a minimal infinite emitter}\}, $$ and the topology has a basis given by the collection $$\{D_{(\beta,B)}: (\beta,B) \in \mathfrak{p}, |\beta|\geq 1\ \} \cup \{D_{(\beta, B),F}:(\beta, B) \in X_{fin}, F\subset \varepsilon\left( B \right), |F|<\infty \},$$ where for each $(\beta,B)\in \mathfrak{p}$ we have that $$D_{(\beta,B)}= \{(\beta, A): A\subset B \text{ and } A\in M_\beta \}\cup\{y \in X: y = \beta \gamma', s(\gamma')\in B\},$$ and, for $(\beta,B)\in X_{fin}$ and $F$ a finite subset of $\varepsilon\left( B \right)$,  $$D_{(\beta, B),F}=  \{(\beta, B)\}\cup\{y \in X: y = \beta \gamma', \gamma_1' \in \ \varepsilon\left( B \right)\setminus F\}.$$
 
\begin{remark}\label{cylindersets} For every $(\beta,B)\in \mathfrak{p}$, we identify $D_{(\beta, B)}$ with $D_{(\beta, B),F}$, where $F=\emptyset$. Furtheremore, we call the basic elements of the topology of $X$ given above by \emph{cylinder sets}.
\end{remark}

It was shown in \cite{GRultrapartial} that the above topological space is metrizable and there a description of convergence of sequences is given. Furthermore, it was shown that, under Condition~(RFUM), the shift space $X$ has a basis of open, compact sets. We recall Condition~(RFUM) below. 


{\bf Condition (RFUM):} For each edge $e\in \GG^1$ its range can be written as $$r(e) = \displaystyle \bigcup_{n=1}^k A_n,$$ where $A_n$ is either a minimal infinite emitter or a single vertex. 

Under Condition (RFUM) we have the following characterization of the range of edges in $\GG$.

\begin{lema}\label{range} Let $\GG$ be an ultragraph with no sinks that satisfies Condition~(RFUM). Then each $A\in \GG^0$ can be written uniquely as $A = \displaystyle \bigcup_{n=1}^k A_n,$ where there exists an unique $k$ such that $|A_k|<\infty$ and $\varepsilon(A_k)< \infty$, $A_j$ is a minimal infinite emitter for $j\neq k$, and $A_j \cap A_k = \emptyset$ for all $j\neq k$. 
\end{lema}
\demo

The existence part follows directly from Condition~(RFUM) and the description of $\GG^0$ given in Lemma~\ref{description}. 

Suppose that $A = \cup_{i=1}^{n} A_i = \cup_{i=1}^m A_i'$, where $|A_k|$ and $|A_k'|<\infty$, $\varepsilon(A_k)$ and $\varepsilon(A_k')< \infty$, and $A_j \cap A_k = \emptyset$ and $A_j' \cap A_k' = \emptyset$ for all $j\neq k$.

Notice that if $A_i \in M_{r(e)}$ and $|A_i|=1$ then $A_i \subseteq A_j'$ for some $j\neq k$. Since $\GG$ has no sinks and $A_j'$ is a minimal infinite emitter, it follows that $A_i =A_j'$. Now, if $A_i \in M_{r(e)}$ and $|A_i|=\infty$, then there exists $j\neq k$ such that $|A_i \cap A_j'|=\infty$. Since there are no sinks, $A_i \cap A_j'$ is a infinite emitter. Then, since both $A_i$ and $A_j'$ are minimal infinite emitters, it follows that $A_i =A_j'$.
Therefore $\cup_{i\neq k} A_i = \cup_{i\neq k} A_i'$ and hence $A_k = A_k'$.

\fim

\begin{corolario}\label{descbasicelemnt} Under the hypothesis of the above lemma, let $(\beta, B)\in \mathfrak{p}$.  
 Then $D_{(\beta,B)}$ can be written as a finite, disjoint union of elements of the collection of cylinder sets $ \{D_{(\beta, B),F}:(\beta, B) \in X_{fin}, F\subset \varepsilon\left( B \right), |F|<\infty \}$, union with the the collection $\{D_{(\alpha,A)}: (\alpha,A) \in \mathfrak{p} \text{ and } \varepsilon(A)<\infty \}$.
\end{corolario}
\demo

Let $(\beta, B)\in \mathfrak{p}$. By Lemma~\ref{range} we have that $B = \displaystyle \bigcup_{n=1}^k A_n,$ where there exists an unique $k$ such that $|A_k|<\infty$ and $\varepsilon(A_k)< \infty$, $A_j$ is a minimal infinite emitter for $j\neq k$, and $A_j \cap A_k = \emptyset$ for all $j\neq k$. 
Let $F_i=\{e\in \GG^1:s(e) \in \cup_{j\neq k} (A_i\cap A_j)\}$ for all $i\neq k$, and let $V = \{s(e): e \in \cup_{i\neq k} F_i\}$. Then
$$ D_{(\beta,B)} = \displaystyle \bigsqcup_{i\neq k} D_{(\beta,A_i),F_i} \bigsqcup_{v\in V} D_{(\beta,\{v\})} \bigsqcup_{v\in A_k} D_{(\beta,\{v\})},$$ and the union is disjoint.

\fim

\begin{remark} It follows from Corollary~\ref{descbasicelemnt} that the collection os sets described in that corollary forms a basis for the topology in $X$.
Notice that if $(\alpha,A) \in \mathfrak{p}$ and $\varepsilon(A)<\infty $, then $D_{(\alpha,A)}= \sqcup_{e \in \varepsilon(A)} D(\alpha e, r(e))$, and the union is disjoint. Furthermore, for $(\alpha,A) \in \mathfrak{p}$ with $\varepsilon(A)<\infty$, and $F\subseteq \varepsilon(A)$, we denote the disjoint union $\sqcup_{e\in \varepsilon(A) \setminus F} D_{(\alpha e,r(e))}$ by $D_{(\alpha,A),F}$.
\end{remark}

Next we prove a result that will be necessary in Section~\ref{kmscrossproduct}, in the construction of a measure associated to a state in $C_0(X)$. 

\begin{proposicao}\label{cilindrosemianel} Let $\GG$ be an ultragraph with no sinks that satisfies Condition~(RFUM). Then the collection of cylinder sets $ \{D_{(\beta, B),F}:(\beta, B) \in X_{fin}, F\subset \varepsilon\left( B \right), |F|<\infty \}$, union with the the collection $\{D_{(\alpha,A)}: (\alpha,A) \in \mathfrak{p} \text{ and } \varepsilon(A)<\infty \}$, union with the empty set, forms a semi-ring.
\end{proposicao}
\demo

Let $S$ denote the collection $ \{D_{(\beta, B),F}:(\beta, B) \in X_{fin}, F\subset \varepsilon\left( B \right), |F|<\infty \}$, union with the the collection $\{D_{(\alpha,A)}: (\alpha,A) \in \mathfrak{p} \text{ and } \varepsilon(A)<\infty \}$, union with the empty set.

Notice that the intersection of two sets in $S$ is again in $S$. 

We have to prove that if $C, C_1 \in S$ are such that $C_1 \subseteq C$, then there exists a finite sequence $C_2, C_3 \ldots C_n \in S$ such that $C$ is equal to the union $\sqcup_{i=1}^n C_i$ and the $C_i$ are disjoint.

Let $C\in S$. Suppose that $C= D_{(\beta, B),F}$, for some $(\beta, B) \in X_{fin}, |F|<\infty$, and let $C_1 \subset C$  be such that $C \setminus C_1 \neq \emptyset$.

Since $C_1 \subseteq C$ we have that $C_1 = D_{(\beta\beta', A),F}$ for some $\beta'$ and $A$ (notice that $\beta' = r(\beta)$ and $F=\emptyset$ are allowed).

Suppose that $C_1=D_{(\beta, A),F'}$, for $(\beta, A)\in X_{fin}$. Then $|\varepsilon(A)|=\infty$ and hence $A=B$. Therefore $C_1=D_{(\beta, B),F'}$. It follows that $$C\setminus C_1 = \displaystyle \bigsqcup_{e\in F'\setminus F} D_{(\beta e, r(e))},$$ 
and, since by Corollary~\ref{descbasicelemnt} each $D_{(\beta e, r(e))}$ is a disjoint union of elements in $S$, we get that $C\setminus C_1$ is a disjoint union of elements in $S$. If $C_1 =D_{(\beta, A)}$, with $|\varepsilon(A)|<\infty$, then $C\setminus C_1 = D_{(\beta, B),F \cup \varepsilon(A)}$. 

Next assume that $C_1 = D_{(\beta\beta', A),F'}$ with $|\beta'|\geq 1$, say $\beta' = \beta_1'\ldots \beta_n'$. Suppose that $(\beta\beta', A)\in X_{fin}$. By Lemma~\ref{range} we can write $r(\beta')=A_1\cup\ldots A_k \ldots A_M$, where $k$ is such that $|A_k|<\infty$ and $\varepsilon(A_k)< \infty$, $A_j$ is a minimal infinite emitter for $j\neq k$, and $A_j \cap A_k = \emptyset$ for all $j\neq k$. 
With this we obtain that 
$$ C\setminus C_1 = D_{(\beta, B),F\cup\{\beta_1'\}}\sqcup D_{(\beta\beta_1', r(\beta_1')),\{\beta_2'\}}\sqcup \ldots \sqcup D_{(\beta\beta_1'\ldots\beta_{n-1}', r(\beta_{n-1}')),\{\beta_n'\}}\bigsqcup$$ 
$$ \bigsqcup_{A_i\neq A} D_{(\beta\beta', A_i),\{\varepsilon((\cup A\cap A_i))\}} \bigsqcup_{e\in F'} D_{(\beta\beta'e,r(e))}.$$
Since by Corollary~\ref{descbasicelemnt} each $D_{(\beta\beta'e,r(e))}$ can be written as a disjoint union of elements in $S$ we get the desired description of $ C\setminus C_1$. 

Now suppose that $C_1 = D_{(\beta\beta', A)}$ with $|\varepsilon(A)|<\infty$. Using the same description of $r(\beta')$ as in the previous case we obtain that 
$$ C\setminus C_1 = D_{(\beta, B),F\cup\{\beta_1'\}}\sqcup D_{(\beta\beta_1', r(\beta_1')),\{\beta_2'\}}\sqcup \ldots \sqcup D_{(\beta\beta_1'\ldots\beta_{n-1}', r(\beta_{n-1}')),\{\beta_n'\}}\bigsqcup$$ 
$$ \bigsqcup_{i\neq k} D_{(\beta\beta', A_i),\{\varepsilon(A)\}} \bigsqcup_{e\in \varepsilon (A_k\setminus A)} D_{(\beta\beta'e,r(e))}.$$

The case when $C= D_{(\beta, B)}$ with $|\varepsilon(B)|<\infty$ is done analogously.

\fim



\subsection{The crossed product construction of an ultragraph C*-algebra}\label{ultracrossproduct}

In this subsection we recall the construction of ultragraph C*-algebras as partial crossed products, as done in \cite{GRultrapartial}. 

Let $\GG$ be an ultragraph with no sinks that satisfy Condition~(RFUM). Denote by $\F$ the free group generated by $\mathcal{G}^1$. We will define a partial action of $\F$ on $X$. For this, let $P\subset \F$ be defined by $$P:=\{e_1...e_n\in \F: e_i\in \mathcal{G}^1: n\geq 1\},$$ 
and define clopen sets $X_c$, for each $c\in \F$, as follows:
\begin{itemize}
\item for the neutral element $0\in \F$ let $X_0=X$;

\item for $a\in P$ define \nl$X_a=\{(\beta,B)\in X_{fin}:\beta_1...\beta_{|a|}=a\}\cup\{\gamma\in \mathfrak{p}^\infty: \gamma_1...\gamma_{|a|}=a\}$; and

$X_{a^{-1}}=\{(A,A)\in X_{fin} :A\subseteq r(a)\}\cup\nl\cup\{(\beta,B)\in X_{fin}: s(\beta)\in r(a)\}\cup \{\gamma\in \mathfrak{p}^\infty:s(\gamma)\in r(a)\};$ 

\item for $a,b\in P$ with $ab^{-1}$ in its reduced form, define \nl $X_{ab^{-1}}=\left\{(a,A)\in X_{fin}:A\subseteq r(a)\cap r(b)\right\}\cup\nl\cup\left\{(\beta,B)\in X_{fin}:\beta_1...\beta_{|a|}=a \text{ and } s(\beta_{|a|+1})\in r(a)\cap r(b)\right\}\cup\nl\cup\left\{\gamma \in \mathfrak{p}^\infty:\gamma_1...\gamma_{|a|}=a \text{ and } s(\gamma_{|a|+1})\in r(a)\cap r(b)\right\}$;

\item for all other $c\in \F$ define $X_c=\emptyset$.

\end{itemize}

\begin{remark} Notice that if $a\in P$ is not a path then $X_a$ is empty. Analogously, if $a,b\in P$ and $r(a)\cap r(b) = \emptyset$ then $X_{ab^{-1}}$ is empty. It is proved in \cite{GRultrapartial} that each $X_c$, with $c\in \F$, is clopen and compact. Furthermore, for each $A\in \mathcal{G}^0$, $X_A:=\{x\in X:s(x)\subseteq A\}$ is clopen and compact.
\end{remark}

Next we recall the definition of the homeomorphisms between the non-empty sets $X_c$, $c\in \F$. For each $a\in P$ such that $X_a$ is non-empty, let $\theta_a:X_{a^{-1}}\rightarrow X_a$ be defined by $\theta_a(x)=a \cdot x$, for each $x\in X_{a^{-1}}$ (here we are using the embedding of $a$ in $\mathfrak{p}$ as $(a,r(a))$). Let $\theta_{a^{-1}}:X_a\rightarrow X_{a^{-1}}$ be defined by 
$\theta_{a^{-1}}((a,A))=(A,A)$, $\theta_{a^{-1}}(ab,B)=(b,B)$ and $\theta(a\gamma)=\gamma$. Finally, for $a,b\in P$ such that $X_{ab^{-1}}$ is non-empty, let $\theta_{ab^{-1}}:X_{ba^{-1}}\rightarrow X_{ab^{-1}}$ be defined by $\theta_{ab^{-1}}(x)=a \cdot \theta_{b^{-1}}(x)$.

\begin{remark} Since, for each $t\in \F$, the map $\theta_t:X_{t^{-1}}\rightarrow X_t$ is a homeomorphism we get that $\alpha_t:C(X_{t^{-1}})\rightarrow C(X_t)$, defined by $\alpha_t(f)=f\circ \theta_{t^{-1}}$, is a *-isomorphism. It follows that $(\{X_t,h_t)\}_{t\in \F}$ is a topological partial action and, consequently, $(\{C(X_t)\}_{t\in \F}, \{\alpha_t\}_{t\in \F})$ is a C*-algebraic partial action of $\F$ in $C(X)$ (see for example \cite{MR3506079}).
\end{remark}


The following result is important for the partial crossed product description of $C^*(\GG)$ and for our work.

\begin{lema}\label{densealgebra} Denote by $1_A$ the characteristic map of $X_A$, for $A\in \GG^0$, and $1_c$ the characteristic map of $X_c$, for $c\in \F$. Then the subalgebra $D\subseteq C_0(X)$ generated by all the characteristic maps $1_c$, $1_A$ and $\alpha_c(1_{c^{-1}}1_A)$, with $c\in \GG^*$ and $A\in\mathcal{G}^0$, is dense in $C_0(X)$. Moreover, for each $0\neq g\in \F$, the subalgebra $D_g$ of $C(X_g)$ generated by all the maps $1_g 1_c$, $1_g 1_A$ and $1_g\alpha_c(1_{c^{-1}}1_A)$, with $c\in \bigcup\limits_{n=1}^\infty (\GG^1)^n$ and $A\in\mathcal{G}^0$, is dense in $C(X_g)$.
\end{lema}

\begin{remark}\label{rmkdensealgebra} Notice that $1_c$ may be written as $\alpha_c(1_{c^{-1}}1_{r(c)})$. Furthermore, $\alpha_c(1_{c^{-1}}1_A)= 1_{r(c)\cap A} \circ \theta_{c^{-1}}$. So the subalgebra $D$ above can be seens as the algebra generated by the characteristic maps $ 1_{r(c)\cap A} \circ \theta_{c^{-1}}$, with $A\in \GG^0$ and $c\in \GG^*$.
\end{remark}

Below we recall the isomorphism between $C^*(\GG)$ and $ C_0(X)\rtimes_\alpha \F$. 

\begin{teorema}\label{crossedproduct} Let $\GG$ be an ultragraph with no sinks that satisfies conditon~(RFUM). Then there exists a bijective *-homomorphism $\Phi:C^*(\mathcal{G})\rightarrow C_0(X)\rtimes_\alpha \F$ such that $\Phi(s_e)=1_e\delta_e$, for each edge $e\in \mathcal{G}^1$, and $\Phi(p_A)=1_A\delta_0$, for all $A\in \GG^0$.
\end{teorema}

We end this section describing the core algebra of $C^*(\G)$ in partial crossed product terms.

\begin{proposicao}\label{isoc0} Let $\GG$ be an ultragraph with no sinks that satisfies Condition~(RFUM). Then $C_0(X)$ is *-isomorphic to the $C^*$-subalgebra of $C^*(\GG)$ generated by $\{s_c p_A s_c^* \mid A\in \GG^0, c\in \GG^* \}$, via a map that takes $1_A$ to $p_A$ and $\alpha_c(1_{c^{-1}}1_A)$ to $s_c p_A s_c^*$, for every $A\in \GG^0$, and every $c\in \GG^*$ with $\ |c|\geq 1$.	
\end{proposicao}
\demo

Recall that $C_0(X)$ is identified with $C_0(X)\delta_0$. By Lemma~\ref{densealgebra}, and the subsequent Remark, the subalgebra $D\subseteq C_0(X)$ generated by all the maps of the form $\alpha_c(1_{c^{-1}}1_A)$, with $c\in \GG^*$ and $A\in\mathcal{G}^0$, is dense in $C_0(X)$. The result now follows since the map $\Phi^{-1}$ of Theorem~\ref{crossedproduct} is a *-isomorphism (over its image) from $C_0(X)\delta_0$ into $C^*(\G)^\gamma$, taking $1_A\delta_0$ to $p_A$ and $\alpha_c(1_{c^{-1}}1_A)\delta_0 = 1_c\delta_c 1_A \delta_A 1_{c^{-1}} \delta_{c^{-1}}$ to $s_c p_A s_c^*$, for every $A\in \GG^0$ and $c\in \bigcup\limits_{n=1}^\infty (\GG^1)^n$.

\fim


\section{KMS states of $C^*$-algebras associated to arbitrary ultragraphs}\label{KMSgeneral}


In this section we describe the set of KMS states of certain one-parameter group of automorphisms of $C^*(\GG)$ in terms of states of the core algebra $C^*(\G)^\gamma$. Throughout this section $\GG$ is an arbitrary ultragraph. 

\subsection*{KMS and ground states.} Suppose that $(A,\R,\rho)$ is a $C^*$-algebraic dynamical system. An element $a$ of $A$ is \emph{analytic} if $t\mapsto \rho_t(a)$ is the restriction of an entire function $z\mapsto \rho_z(a)$ on $\C$. A state $\phi$ of $(A,\R,\rho)$ is a \emph{KMS state with inverse temperature $\beta$} (or a KMS$_\beta$ state) if $\phi(ab)=\phi(b\rho_{i\beta}(a))$ for all analytic elements $a,b$.  
It is usually enough  to check the KMS condition  on a set of analytic elements which span a dense subalgebra of $A$. A state $\phi$ on $A$ is a \emph{ground state} of $(A,\R,\rho)$ if for every $a, b$ analytic in $A$, the entire function $z\mapsto \phi(a \rho_z(b))$ is bounded on the upper-half plane.

\subsection*{Generalized gauge action.}\label{Generalgaugeaction}
Let $\G$ be an ultragraph and let  $N$ be a positive function on $\mathcal{G}^1$ such that there is a constant $K$ such that $N(e)>K$ for all $e\in \G^1$.
Extend the function $N$ to $N:\mathcal{G}^*\rightarrow\R_+^*$ by defining $N(A)=1$, for $A\in \G^0$, and $N(\lambda)=N(e_1)\dots N(e_m)$ for $\lambda=e_1\dots e_m \in \mathcal{G}^*$, with $|\lambda|>0$. 

The following lemma shows that the function $N$ gives an action of $\R$ on $C^*(\G)$.

\begin{lema}\label{action}
 Let $\G$ be an ultragraph and let  $N:\mathcal{G}^1\rightarrow\R_+^*$ as above. Then there is a strongly continuous  action
 $\rho^c:\R\rightarrow \Aut(C^*(\G))$  such that $\rho^c_t(p_A)=p_A$, for all $A\in \G^0$, and $\rho^c_t(s_e)=N(e)^{it}s_e$ for $e\in \G^1$.
\end{lema}
\demo
Fix $t\in \R$.  For each $e\in \G^1$, define $T_e:=N(e)^{it}s_e$.
Since $T_eT_e^*=s_es_e^*$ and $T_e^*T_e=s_e^*s_e$, it is easy to see that $\{T_e:e\in \G^1\}$ and $\{p_A:A\in \G^0\}$ satisfies the conditions of Definition~\ref{def of C^*(mathcal{G})}. Now the universal property gives a homomorphism $\rho_t^c:C^*(\mathcal{G})\rightarrow C^*(\mathcal{G})$ such that $\rho^c_t(p_A)=p_A$, for all $A\in \G^0$, and $\rho^c_t(s_e)=N(e)^{it}s_e$ for $e\in \G^1$.
Observe that  $\rho^c_{t}\circ \rho^c_{t'}=\rho^c_{t+t'}$ for $t,t'\in \R$. Also notice  that the identity map on $C^*(\G)$ is $\rho^c_0$. It follows that $(\rho_{t}^c)^{-1}=\rho_{-t}^c$ and hence $\rho_{t}^c\in \Aut(C^*(\G))$.
	Therefore  $\rho^c$ is a homomorphism of $\R$ into the $\Aut(C^*(\G))$. 
	
	To see that $\rho^c$ is strongly continuous  we must prove that  $t\mapsto \rho_t^c(b)$ is continuous for all $b\in C^*(\G)$. Fix $\epsilon>0$ and $b\in C^*(\G)$. There is a linear combination $d$ of generators in $C^*(\G)$ such that $\|b-d\|<\frac{\epsilon}{3}$. Since $t\mapsto \rho_t^c(d)$ is continuous, there exists some $\delta>0$
	such that $|t-u|<\delta \Rightarrow \|\rho_t^c(d)-\rho_u^c(d)\|<\frac{\epsilon}{3}$.
	Now for $|t-u|<\delta$ we have that
	\[\|\rho_u^c(b)-\rho_t^c(b)\|\leq \|\rho_u^c(b-d)\|+\|\rho_u^c(d)-\rho_t^c(d)\|+\|\rho_t^c(b-d)\|<\epsilon,\]
	as required.
\fim

\subsection*{Some properties of $C^*(\G)$.}

Given an ultragraph $\G$,   it is shown in  \cite[page 351]{MR2050134} that 
\[C^*(\G)=\lsp \{s_\mu p_A s_\nu^*:\mu,\nu\in \G^* , A\in \G^0\}.\]
Now let  $N:\mathcal{G}^1\rightarrow\R_+^*$ and let  $\rho$ be the associated action as in Lemma~\ref{action}. For every $s_\mu p_A s_\nu^*\in C^*(\G)$, the map 
$t\mapsto\rho_t(s_\mu p_A s_\nu^*)=N(\mu)^{it}N(\nu)^{-it}s_\mu p_A s_\nu^*$ on $\R$ extends to an analytic function on all the complex plane. Thus there are analytic elements which span a dense subalgebra of  $C^*(\G)$, and therefore to study KMS states one can consider only these elements.  

Also it is shown in \cite[page 350]{MR2050134}  that for $\mu,\nu \in \G^*$, with $|\mu|, |\nu| \geq 1,$ we have

\begin{align}\label{productformula}
	s_\nu^*s_\mu&=\begin{cases}
	s_{\nu'}^* \notag&\text{if $\nu=\mu\nu',\nu'\notin\G^0$}\\
	s_{\mu'} \notag&\text{if $\mu=\nu\mu',\mu'\notin\G^0$}\\
	p_{ r(\nu)} \notag&\text{if $\mu=\nu$}\\
	0\notag&\text{otherwise.}
	\end{cases}\\
\end{align}

There is also a strongly continuous gauge action $\gamma:\T\rightarrow \Aut(C^*(\G))$	such that $\gamma_z(s_\mu)=zs_\mu$ and $\gamma_z(p_A)=p_A$. 
The  \textit{core algebra} $C^*(\G)^\gamma$ is the fixed point subalgebra  for the gauge action.

The proof of the next two lemmas works the same as for graph C*-algebras as done in Raeburn's book \cite{MR2135030}.

\begin{lema}\label{fixed-alg}	
\[C^*(\G)^\gamma=\lsp \{s_\mu p_A s_\nu^*:\mu,\nu\in \G^* , A\in \G^0,|\mu|=|\nu|\}.\]
\end{lema}

\begin{lema}\label{cond-exp}	
There is a conditional expectation $\Psi:C^*(\G)\rightarrow C^*(\G)^\gamma$ such that 
\[\Psi(s_\mu p_A s_\nu^*)=\delta_{|\mu|,|\nu|} s_\mu p_A s_\nu^*,\] for all $\mu,\nu\in \G^* , A\in \G^0$. 
\end{lema}

\subsection*{KMS states for the generalized gauge action}	
\begin{lema}
Let $N:\mathcal{G}^1\rightarrow\R_+^*$ be such that $N(\mu)\neq 1$ for $\mu\in \G^*\setminus\G^0$ and let  $\rho$ be the action of Lemma~\ref{action}. Let
  $\beta\in \R $. Suppose that $\phi, \phi'$ are KMS$_\beta$ states on $C^*(\G)$ coinciding on the core algebra $C^*(\G)^\gamma$. Then $\phi=\phi'$. 
\end{lema}
\demo
Take an  element $s_\mu p_As_\nu^*\in C^*(\G)$ with $|\mu|\neq|\nu|$ and let $\omega:=\phi-\phi'$. We aim to show that $\omega(s_\mu p_As_\nu)=0$. Since $\omega$ is a KMS$_\beta$ state, by applying the KMS condition, we have that
\begin{align*}
	\omega(s_\mu p_As_\nu^*)&=\omega(p_As_\nu^*\rho_{i\beta}(s_\mu) )=N(\mu)^{-\beta}\omega(p_As_\nu^* s_\mu )\\
	&=\begin{cases}
   N(\mu)^{-\beta}p_As_{\mu'} &\text{if $\mu=\nu\mu',\mu'\notin\G^0$}\\
	 N(\mu)^{-\beta}p_As_{\nu'}^* &\text{if $\nu=\mu\nu',\nu'\notin\G^0$}\\
	 0&\text{otherwise.}\\
	\end{cases}
\end{align*}
Now it suffices to prove that $\omega(p_Bs_\lambda)=\omega(p_Bs_\lambda^*)=0$ for all $B\in \G^0, \lambda\in \G^*\setminus\G^0$. To prove this, notice that if $C^*(\G)$ has a unit, then
\[\omega(p_Bs_\lambda)=\omega(p_Bs_\lambda 1)=N(\lambda)^{-\beta}\omega( 1 p_Bs_\lambda)=N(\lambda)^{-\beta}\omega(p_Bs_\lambda).\]
Since $N(\lambda)\neq 1$, it follows that $\omega(p_Bs_\lambda)=0$ as required.

If there is no unit for $C^*(\G)$, we apply the same argument for the approximate unit $\{p_A: A\in \GG^0\}$.

\fim

\begin{proposicao}\label{equivalenciaAE}
Let $N:\mathcal{G}^1\rightarrow\R_+^*$ be such that $N(\mu)\neq 1$ for $\mu\in \G^*\setminus\G^0$ and let  $\rho$ be the action of Lemma~\ref{action}. Suppose that
  $\beta\in \R $ and  $\phi$ is a state on $C^*(\G)$. Then the restriction $\psi:=\phi|_{C^*(\G)^\gamma}$ satisfies 
\begin{equation}\label{kms}
\psi(s_\mu p_As_\nu^*)=\delta_{\mu,\nu}N(\mu)^{-\beta}\psi(p_{A\cap r(\mu)}).
\end{equation}
Conversely, for any state  $\psi$  on ${C^*(\G)^\gamma}$ satisfying \eqref{kms}, $\phi=\psi\circ \Psi$ is a KMS$_\beta$ state on $C^*(\G)$, where $\Psi$ is the conditional expectation as in Lemma~\ref{cond-exp}. Furthermore the obtained correspondence is an affine bijection. 
\end{proposicao}
\demo
Suppose that $\phi$ is a KMS$_\beta$ state on $C^*(\G)$ and let $\psi$ be its restriction to ${C^*(\G)^\gamma}$. Take $s_\mu p_As_\nu\in C^*(\G)$ with $|\mu|\neq|\nu|$. Then the KMS condition impies that
\begin{align*}
	\phi(s_\mu p_As_\nu^*)=\phi(p_As_\nu^* \rho_{i\beta}(s_\mu))= N(\mu)^{-\beta}\phi(p_As_\nu^* s_\mu ).
\end{align*}
Therefore,  by formula~\eqref{productformula}, we have that $\psi(s_\mu p_As_\nu^*)=\delta_{\mu,\nu}N(\mu)^{-\beta}\psi(p_{A\cap r(\mu)}).$

Next suppose that $\psi$ is a state on $C^*(\G)^\gamma$ which satisfies \eqref{kms}.  To see that $\phi=\psi\circ\Psi$ is a KMS$_\beta$ state, it suffices to prove the KMS condition
\begin{equation}\label{kmscondition}
\phi(ab)=N(\mu)^{-\beta}N(\nu)^{\beta}\phi(ba),
\end{equation}
for $a=s_\mu p_As_\nu^*$, $b=s_\lambda p_B s_\tau^*$ where $\mu,\nu,\lambda,\tau\in \G^*$ and $A,B\in \G^0$. To see this first note that 
\begin{align*}
	ab=s_\mu p_As_\nu^*s_\lambda p_B s_\tau^*&=\begin{cases}
	s_\mu p_As_{\nu'}^* p_B s_\tau^*&\text{if $\nu=\lambda\nu',\nu'\notin\G^0$}\\
	s_\mu p_As_{\lambda'} p_B s_\tau^*&\text{if $\lambda=\nu\lambda',\lambda'\notin\G^0$}\\
	s_\mu p_{A\cap r(\nu)\cap B} s_\tau^*&\text{if $\lambda=\nu$}\\
	0&\text{otherwise.}\\
	\end{cases}\\
	&=\begin{cases}
	s_\mu p_As_{\nu'}^*s_\tau^*&\text{if $\nu=\lambda\nu',\nu'\notin\G^0,s(\nu')\in B$}\\
	s_\mu s_{\lambda'}p_B  s_\tau^*&\text{if $\lambda=\nu\lambda',\lambda'\notin\G^0, s(\lambda')\in A$}\\
	s_\mu p_{A\cap r(\nu)\cap B} s_\tau^*&\text{if $\lambda=\nu$}\\
	0&\text{otherwise.}\\
	\end{cases}\\
	&=\begin{cases}
	s_\mu p_As_{\tau\nu'}*&\text{if $\nu=\lambda\nu',\nu'\notin\G^0,s(\nu')\in B$}\\
	s_{\mu \lambda'} p_B s_\tau^*&\text{if $\lambda=\nu\lambda',\lambda'\notin\G^0, s(\lambda')\in A$}\\
	s_\mu p_{A\cap r(\nu)\cap B} s_\tau^*&\text{if $\lambda=\nu$}\\
	0&\text{otherwise.}\\
	\end{cases}\\
\end{align*}

Similarly 

\begin{align*}
	ba=s_\lambda p_B s_\tau^*s_\mu p_As_\nu^*
	&=\begin{cases}
	s_\lambda p_B s_{\nu\tau'}^*&\text{if $\tau=\mu\tau',\tau'\notin\G^0, s(\tau')\in A$}\\
	s_{\lambda\mu'} p_As_\nu^*&\text{if $\mu=\tau\mu',\mu'\notin\G^0, s(\mu')\in B$}\\
	s_\lambda p_{B\cap r(\mu)\cap A} s_\nu^*&\text{if $\tau=\mu$}\\
	0&\text{otherwise.}\\
	\end{cases}\\
\end{align*}

Now we apply  $\phi=\psi\circ\Psi$ to both $ab$ and $ba$. Then 
\begin{align*}
	\phi(ab)&=\begin{cases}
	N(\mu)^{-\beta}\psi(p_{A\cap r(\mu)})&\text{if $\nu=\lambda\nu',\nu'\notin\G^0,s(\nu')\in B,\mu=\tau\nu'$}\\
	 N(\tau)^{-\beta}\psi(p_{B\cap r(\tau)}) &\text{if $\lambda=\nu\lambda',\lambda'\notin\G^0, s(\lambda')\in A,\mu \lambda'=\tau $}\\
	N(\mu)^{-\beta}\psi(p_{A\cap r(\nu)\cap B\cap r(\mu)}) &\text{if $\lambda=\nu, \mu=\tau$}\\
	0&\text{otherwise,}\\
	\end{cases}\\
\end{align*}
\begin{align*}
	\phi(ba)&=\begin{cases}
	N(\lambda)^{-\beta}\psi(p_{B\cap r(\lambda)})&\text{if $\tau=\mu\tau',\tau'\notin\G^0, s(\tau')\in A,\lambda=\nu\tau'$}\\
	N(\nu)^{-\beta} \psi(p_{A\cap r(\nu)})&\text{if $\mu=\tau\mu',\mu'\notin\G^0, s(\mu')\in B,\lambda\mu'=\nu$}\\
	N(\nu)^{-\beta} \psi(p_{B\cap r(\mu)\cap A\cap r(\nu)})&\text{if $\tau=\mu,\lambda=\nu$}\\
	0&\text{otherwise.}\\
	\end{cases}\\
\end{align*}
Now if $\nu=\lambda\nu',s(\nu')\in B,\mu=\tau\nu'$, then $r(\mu)=r(\nu')=r(\nu)$ and hence
\[N(\mu)^{-\beta}N(\nu)^{\beta}\phi(ba)=N(\mu)^{-\beta}\psi(p_{A\cap r(\nu)})=N(\mu)^{-\beta} \psi(p_{A\cap r(\mu)})=\phi(ab).\]
If $\lambda=\nu\lambda',s(\lambda')\in A,\mu \lambda'=\tau $, then $r(\tau)=r(\lambda')=r(\lambda)$. Also we have $\frac{N(\lambda)}{N(\nu)}=N(\lambda')=\frac{N(\tau)}{N(\mu)}$ and hence $N(\lambda)^{-\beta}N(\mu)^{-\beta}N(\nu)^{\beta}=N(\tau)^{-\beta}$. Then 
\[N(\mu)^{-\beta}N(\nu)^{\beta}\phi(ba)=N(\tau)^{-\beta}\psi(p_{B\cap r(\lambda)})=N(\tau)^{-\beta} \psi(p_{A\cap r(\tau)})=\phi(ab).\]
Finally if $\lambda=\nu, \mu=\tau$, then clearly  
\[N(\mu)^{-\beta}N(\nu)^{\beta}\phi(ba)=\phi(ab),\]
and  we have proven \eqref{kmscondition}.
\fim


\section{KMS states of ultragraph C*-algebras realized as partial crossed products}\label{kmscrossproduct}


In this section we will further describe the set of KMS states associated to the one-parameter group of automorphisms described in the previous section. To this end we will make use of the description of C*-algebras associated to ultragraphs that have no sinks and satisfy Condition~(RFUM) as partial crossed products (see Section~\ref{ultracrossproduct}) and build from the ideas in Section~4 of \cite{MR3539347}. Before we proceed we set up some notation and recall the construction of the one-parameter group of automorphisms via partial crossed product theory. 

{\bf Assumption:} From now on all ultragraphs are assumed to have no sinks.

%

Recall from \cite[Theorem 4.3]{MR1953065} that given any function $N:\GG^1\to (1,\infty)$ there exists a unique strongly continuous one-parameter group $\sigma$ of automorphisms of $C_0(X)\rtimes \F$ such that
\begin{equation}\label{eq:sigma-from-N}
  \sigma_t(b)=\bigl(N(e)\bigr)^{it}b\,\text{ and }\,\sigma_t(c)=c
\end{equation}
for all $e\in \GG^1$, all $b\in C(X_e) \delta_e$, and all $c\in C_0(X)\delta_0$.

If  $N(e)=\exp(1)$ for every $e\in E^1$, then $\sigma_t$ is
$2\pi$-periodic, and so induces a strongly continuous action $\beta:\T\to\aut(C_0(X)\rtimes \F)$
such that $\beta_z(1_e\delta_e)=z1_e \delta_e\text{ and }\beta_z(f\delta_0)=f\delta_0$
for all $z\in\T$, $e\in \GG^1$, and $f\in C_0(X)$. Alternatively, one can build such an action proceeding as in the proof of \cite[Theorem 4.11]{GRultrapartial}.

We then have the following

\begin{corolario} Let $\GG$ be an ultragraph with no sinks that satisfies condition (RFUM), let $\gamma$ be the gauge action on $C^*(\GG)$, $\beta$ as above, and let $\Phi:C^*(\mathcal{G})\rightarrow C_0(X)\rtimes_\alpha \F$ be the isomorphism of Theorem \ref{crossedproduct}. Then $$\Phi \circ\gamma_z=\beta_z\circ\Phi$$ for all $z\in\T$.
\end{corolario}

Proceeding as in \cite{MR3539347}, given a function $N:\GG^1\to (1,\infty)$, and letting $\sigma$ be the unique strongly continuous one-parameter group of automorphisms of $\C(X)\times \F$ given by \eqref{eq:sigma-from-N} we have, by the Corollary above and Theorem \ref{crossedproduct}, a unique strongly continuous one-parameter group $\sigma$ of automorphisms of $C^*(\GG)$ such that
\begin{equation*}
  \sigma_t(s_e)=\bigl(N(e)\bigr)^{it}s_e \text{ and }\sigma_t(p_A)=p_A
\end{equation*}
for all $e\in \GG^1$ and $A\in \GG^0$.

\begin{remark} In Section~\ref{Generalgaugeaction} we considered any positive function $N:\GG^1\to \R^+_*$ such that $N(e)>K$ for some $K>0$. In this section we are restricting our attention to functions $N:\GG^1\to (1,\infty)$. Furthermore, the one-parameter group of automorphisms considered above is of the same form as the one-parameter group of
automorphisms considered in Lemma \ref{action}. 
\end{remark}

Next we set up some notation that will simplify our follow up statements. 

For $0\leq\beta<\infty$ we define the following sets:
\begin{itemize}
		\item[$A^\beta$:] the set of KMS$_\beta$ states for $C^*(\GG)$,
		\item[$B^\beta$:] the set of states $\omega$ of $C_0(X)$ that satisfy the scaling condition $\omega(f\circ\theta_e^{-1})=N(e)^{-\beta}\omega(f)$ for all $e\in \GG^1$ and all $f\in C_0(X_{e^{-1}})$,
	  \item[$C^\beta$:] the set of regular Borel probability measures $\mu$ on $X$ that satisfy the scaling condition $\mu(\theta_e(A))=N(e)^{-\beta}\mu(A)$ for every $e\in \GG^1$ and every Borel measurable subset $A$ of $X_{e^{-1}}$, 
	  \item[$D^\beta$:] the set of
functions $m:\GG^0\to [0,1]$ satisfying 
\begin{enumerate}\renewcommand{\theenumi}{m\arabic{enumi}}
  	\item\label{item:n1} $\lim_{A\in \GG^0} m(A) =1$;
    \item\label{item:n2} $m(A)=\sum_{e: s(e)\in A}N(e)^{-\beta}m\bigl(r(e)\bigr)$ if $|\varepsilon(A)|<\infty$;
    \item\label{item:n3} $m(A)\geq \sum_{e\in F}N(e)^{-\beta}m\bigl(r(e)\bigr)$ for every finite subset $F$ of $\varepsilon(A)$;
    \item\label{item:n4} $m(A\cup B) = m(A)+m(B)-m(A\cap B)$, and
  \end{enumerate}
  \item[$E^\beta$:] the set of states $\psi$  on ${C^*(\G)^\gamma}$ satisfying $
\psi(s_c p_A s_d^*)=\delta_{c,d}N(c)^{-\beta}\psi(p_{A\cap r(c)})$, for all $c,d \in \G^*$ and $A \in \GG^0$.
	\end{itemize}

Our goal for the reminder of this section is to prove the following theorem.

\begin{teorema}\label{teoremaprincipal} Let $\G$ be an ultragraph with no sinks that satisfy Condition~(RFUM). Then there exists a convex isomorphism between $A^\beta$, $B^\beta$, $C^\beta$, $D^\beta$, and $E^\beta$.

\end{teorema}

By Proposition~\ref{equivalenciaAE} we already have a convex isomorphism between $A^\beta$ and $E^\beta$. To prove the other equivalences we need to prove a few auxiliary results first.

\begin{proposicao}\label{From B to D}
Let $\GG$ be an ultragraph with no sinks that satisfies Condition~(RFUM) and let $M$ be a function from $\GG^1$ to $[0,1]$. Then there is a convex, injective map between the set of states $\omega$ of $C_0(X)$ such that $\omega(f\circ\theta_e^{-1})=M(e)\omega(f)$ for all $e\in \GG^1$ and all $f\in C_0(X_{e^{-1}})$, and the set of
functions $m:\GG^0\to [0,1]$ satisfying 
\begin{enumerate}\renewcommand{\theenumi}{m\arabic{enumi}'}
  	\item\label{item:n1linha} $\lim_{A\in \GG^0} m(A) =1$;
    \item\label{item:n2linha} $m(A)=\sum_{e: s(e)\in A}M(e)m\bigl(r(e)\bigr)$ if $|\varepsilon(A)|<\infty$;
    \item\label{item:n3linha} $m(A)\geq \sum_{e\in F}M(e)m\bigl(r(e)\bigr)$ for every finite subset $F$ of $\varepsilon(A)$;
    \item\label{item:n4linha} $m(A\cup B) = m(A)+m(B)-m(A\cap B)$.
  \end{enumerate}
Furthermore, the correspondence takes a state $\omega$ to the function defined, for $A\in \GG^0$, by $m(A)= \omega(1_A)$.
 
\end{proposicao}

\demo

Let $\omega$ be a state of $C_0(X)$ such that $\omega(f\circ\theta_e^{-1})=M(e)\omega(f)$ for all $e\in \GG^1$ and all $f\in C_0(X_{e^{-1}})$. Let $m$ be the function from $\GG^0$ to $[0,1]$ given by
  \begin{equation*}
  	m(A)=\omega(1_A).
  \end{equation*}

First we prove \eqref{item:n1linha}: Notice that $\{A:A\in \GG^0 \}$ is a directed set (with inclusion as preorder), and $\{1_A:A\in \GG^0 \}$ is an increasing approximate unity for $C_0(X)$. Then, by Theorem 3.3.3 of \cite{MR1074574}, we have that $1= \lim \omega(1_A) = \lim m(A)$ and hence $m$ satisfies \eqref{item:n1}.

Before we show \eqref{item:n2linha} and \eqref{item:n3linha} notice that if $e\in \GG^1$, then
$$    \omega\bigl(1_e \bigr)=
    \omega\bigl(1_{e^{-1}}\circ \theta_e^{-1} \bigr)= M(e) \omega (1_{e^{-1}}) = M(e) \omega \bigl(1_{r(e)}\bigr) = M(e) m \bigl(r(e)\bigr).$$

Proof of \eqref{item:n2linha}: Notice that if $|\varepsilon(A)|<\infty$ then $1_A = \sum_{e: s(e)\in A} 1_e$. Hence 
  \begin{equation*}
      m(A)=\omega\bigl(1_A \bigr)
      =\sum_{e: s(e)\in A}\omega\bigl(1_e\bigr)
      =\sum_{e: s(e)\in A} M(e) m \bigl(r(e)\bigr),
  \end{equation*}
which gives \eqref{item:n2linha}.  

Proof of \eqref{item:n3linha}: Suppose that $F$ is a finite subset $\varepsilon(A)$. Then $1_A \geq \sum_{e\in F} 1_e$ and we have that
  \begin{equation*}
      m(A)=\omega\bigl(1_A \bigr) \geq \sum_{e\in F}\omega\bigl(1_e\bigr)
      =\sum_{e\in F}M(e)m\bigl(r(e)\bigr).
  \end{equation*}

Proof of \eqref{item:n4linha}: This follows from the linearity of $\omega$.

Next we prove that the correspondence given above is injective. 

Let $\omega_1$ and $\omega_2$ be states such that $\omega_i(f\circ\theta_e^{-1})=M(e)\omega_i(f)$ for all $e\in \GG^1$ and all $f\in C_0(X_{e^{-1}})$, $i=1,2$. Suppose that $\omega_1(1_A) = \omega_2(1_A)$ for all $A\in \GG^0$. We have to show that $\omega_1=\omega_2$. By Remark~\ref{rmkdensealgebra} it is enough to show that $\omega_1\bigl(1_{r(c)\cap A} \circ \theta_{c^{-1}}\bigr)= \omega_2\bigl(1_{r(c)\cap A} \circ \theta_{c^{-1}}\bigr)$ for all $A\in \GG^0$ and $c\in \GG^*$, $|c|\geq 1$. This follows from the following computations for $c=c_1\ldots c_n$ and $A \in \GG^0$:

\begin{align*}
   \omega_1\bigl( 1_{r(c)\cap A} \circ \theta_{c_n^{-1}\ldots c_1^{-1}} \bigr)
    & = \omega_1( 1_{\theta_{c_n} (r(c)\cap A)} \circ \theta_{c_{n-1}^{-1}\ldots c_1^{-1}} \bigr)
    = \ldots \\
	& =	\omega_1\bigl( 1_{\theta_{c_2\ldots c_n} (r(c)\cap A)} \circ \theta_{c_1^{-1}}\bigr)
 = M(c_1) \omega_1 \bigl(  1_{\theta_{c_2\ldots c_n}} \bigr) \\ 	
 & = M(c_1) \omega_1 \bigl(  1_{\theta_{c_3\ldots c_n}} \circ \theta_{c_2}^{-1} \bigr)=\ldots	\\
    & = M(c_1)\ldots M(c_n) \omega_1 \bigl(1_{r(c)\cap A} \bigr) \\
        & = M(c_1)\ldots M(c_n) \omega_2 \bigl(1_{r(c)\cap A} \bigr) = \ldots \\
        & = \omega_2\bigl( 1_{r(c)\cap A} \circ \theta_{c_n^{-1}\ldots c_1^{-1}} \bigr).
  \end{align*}

\fim

Next we will build a measure on the shift space $X$ following some of the ideas in section 5.5 of \cite{MR2356043}. Recall that, as mentioned in Remark~\ref{cylindersets}, we are writing $D_{(\alpha,A),\emptyset}$ for a cylinder of the form $D_{(\alpha,A)}$, so that all generalized cylinders can be written as $D_{(\beta,B),F}$. 

Given a function $M$ from $\GG^1$ to $[0,1]$ and a function $m:\GG^0\to [0,1]$ satisfying \ref{item:n1linha} to \ref{item:n4linha} of Proposition~\ref{From B to D}, we extend $M$ to $\G^*$ to the function $M:\G^* \rightarrow [0,1]$ given by $M(A)= m(A)$, for $A\in \G^0$, and $M(\beta)=M(e_1)\ldots M(e_n)$ for $e_1...e_n \in \G^*$, with $|\beta|\geq 1$.  Now define a function on the set of all generalized cylinders by
\begin{equation}\label{defkappa} \kappa(D_{(\beta,B),F})=M(\beta)m(B)-\sum_{e\in F}M(\beta e)m(r(e)) 
\end{equation}
which, by m3', takes values on non-negative numbers. We also define $\kappa(\emptyset)=0$.

\begin{lema}\label{measureonsemiring}
Let $\G$ be an ultragraph with no sinks that satisfy Condition~(RFUM) and $\kappa$ be the function defined by Equation~(\ref{defkappa}). Then the restriction of $\kappa$ to the semi-ring $S$ given by Proposition~\ref{cilindrosemianel} is a measure such that $\kappa(\theta_e(V))=M(e)\kappa(V)$, for every $e\in \GG^1$ and every subset $V$ of $X_{e^{-1}}\cap S$.
\end{lema}
\demo

We have to prove that $\kappa$ is countably additive on $S$. Since all elements of $S$ are compact open sets, it is actually sufficient to show that $\kappa$ is additive on $S$. Suppose that 

\begin{equation}\label{eq.disjoint.union.cylinder}
D_{(\beta,B),F}=\bigsqcup_{i=1}^n D_{(\beta_i,B_i),F_i}.
\end{equation}

Notice that in this case, for all $i$, $\beta_i=\beta\beta'_i$ for some $\beta'_i$. We use an induction argument on $m=\max_i\{|\beta_i|-|\beta|\}$.

First suppose that $m=0$ so that $\beta_i=\beta$ and $B_i\subseteq B$ for all $i$. In the case that $(\beta,B)\in \mathfrak{p}$ with $|\varepsilon(\beta)|<\infty$, by our convention $F=\emptyset$. Also, since $B_i\subseteq B$, then  $|\varepsilon(\beta_i)|<\infty$ and $F_i=\emptyset$ for all $i$. It follows that $B=\sqcup_i B_i$. By m4'
\begin{align*}
\kappa(D_{(\beta,B)})&=M(\beta)m(B)  \\
&=M(\beta)m\left(\sqcup_i B_i\right) \\
&=M(\beta)m\sum_i m(B_i) \\
&=\sum_i \kappa(D_{(\beta_i,B_i)}).
\end{align*}

For the case that $(\beta,B)\in X_{fin}$, since $(\beta,B)\in D_{(\beta,B)}$, we must have $B_{i_0}=B$ for some $i_0$, and so $F_{i_0}\supseteq F$. For simplicity, suppose that $i_0=1$. Notice that for $i\neq 1$, $B_i$ cannot be a infinite emitter so that $F_i=\emptyset$. In order for (\ref{eq.disjoint.union.cylinder}) to hold, it is necessary that 
$\sqcup_{i\geq 2}\varepsilon(B_i) = F_1\setminus F$. Using the definition of $\kappa$ and m2' we obtain that
\begin{align*}
\kappa(D_{(\beta,B)},F)&=M(\beta)m(B)-\sum_{e\in F}M(\beta e)m(r(e))  \\
&=M(\beta)m(B)-\sum_{e\in F_1}M(\beta e)m(r(e))+\sum_{e\in F_1 \setminus F}M(\beta e)m(r(e)) \\
&=M(\beta)m(B)-\sum_{e\in F_1}M(\beta e)m(r(e))+\sum_{i=2}^{n}\sum_{e\in\varepsilon(B_i)}M(\beta e)m(r(e)) \\
&=M(\beta)m(B)-\sum_{e\in F_1}M(\beta e)m(r(e))+\sum_{i=2}^{n}M(\beta)m(B_i) \\
&=\sum_i \kappa(D_{(\beta_i,B_i),F_i}).
\end{align*}

Now suppose that the result holds for all $0\leq k <m$. Let us do the case where $(\beta,B)\in X_{fin}$, the other being analogous (CHECK). As above, we must have that there exists $i_0$ such that $\beta_{i_0}=\beta$, $B_{i_0}=B$ and $F_{i_0}\supseteq F$. Again, for simplicity, suppose that $i_0=1$. Also suppose that for some $i_1$, $\beta_i=\beta$ for $2\leq i \leq i_1$ and $\beta_i\neq \beta$ for $i_1 < i \leq n$. As in the case $m=0$, $F_i=\emptyset$ for $2\leq i \leq i_1$, but now $\sqcup_{2\leq i \leq i_1}\varepsilon(B_i) \subseteq F_1\setminus F$. Define $C=s((F_1\setminus F)\setminus \sqcup_{2\leq i \leq i_1}\varepsilon(B_i))$ so that $\varepsilon(C)$ is finite. Notice that
\[D_{(\beta,B),F}=D_{(\beta,C)}\sqcup\bigsqcup_{i=1}^{i_1} D_{(\beta_i,B_i),F_i}\]
and
\[D_{(\beta,C)}=\bigsqcup_{i=i_1+1}^{n} D_{(\beta_i,B_i),F_i}.\]
Also, by the case $m=0$, we have that 
\begin{equation}\label{eq:sum0}
\kappa(D_{(\beta,B),F})=\kappa(D_{(\beta,C)})+\sum_{i=1}^{i_1} \kappa (D_{(\beta_i,B_i),F_i}).
\end{equation}

Since $\varepsilon(C)$ is finite then $D_{(\beta,C)}=\bigsqcup_{s(e)\in C}D_{(\beta e,r(e))}$. Using Corollary~\ref{descbasicelemnt}, m2' and m4', we can write
\[D_{(\beta,C)}=\bigsqcup_j D_{(\beta e_j, A_j)},\]
where $D_{(\beta e_j, A_j)}\in S$, in such way that
\begin{equation} \label{eq:sum1}\kappa(D_{(\beta,C)})=\sum_j\kappa( D_{(\beta e_j, A_j)}).\end{equation}

Now
\[D_{(\beta e_j, A_j)}=D_{(\beta e_j, A_j)}\cap D_{(\beta,C)}=\bigsqcup_{i=i_1+1}^{n}D_{(\beta e_j, A_j)}\cap D_{(\beta_i,B_i),F_i},\]
and notice that $D_{(\beta e_j, A_j)}\cap D_{(\beta_i,B_i),F_i} = D_{\gamma_j,D_j}$ for some element of $S$ such that $|\gamma_j|=|\beta_i|$ (since $|\beta_i|\geq |\beta|+1=|\beta e_j|$). This implies that $\max_j\{|\gamma_j|-|\beta e_i|\}<m$. By the induction hypothesis
\begin{equation} \label{eq:sum2}\kappa(D_{(\beta e_j, A_j)})=\sum_{i=i_1+1}^{n}\kappa(D_{(\beta e_j, A_j)}\cap D_{(\beta_i,B_i),F_i}).\end{equation}

On the other hand, for $i_1+1\leq i\leq n$, 
\[D_{(\beta_i,B_i),F_i}=D_{(\beta_i,B_i),F_i}\cap D_{\beta,C}=\bigsqcup_j D_{(\beta e_j, A_j)}\cap D_{(\beta_i,B_i),F_i},\]
and using the case $m=0$ (even if some intersections are empty) we obtain that
\begin{equation} \label{eq:sum3}\kappa(D_{(\beta_i,B_i),F_i})=\sum_j \kappa(D_{(\beta e_j, A_j)}\cap D_{(\beta_i,B_i),F_i}). \end{equation}

Putting equations (\ref{eq:sum1}), (\ref{eq:sum2}) and (\ref{eq:sum3}) together we arrive at
\[\kappa(D_{(\beta,C)})=\sum_{i=i_1+1}^{n} \kappa(D_{(\beta_i,B_i),F_i}),\]
and, from (\ref{eq:sum0}), we get that $\kappa$ is a measure.

Finally, to prove the scaling condition, notice that if $V=D_{(\beta,B),F}$ is an element of $S$ contained in $X_{e^{-1}}$ then $\theta_e(V)=D_{(e\beta,B),F}$. Hence, from the definition of $\kappa$, and the multiplicativity of $M$, the result follows. 
\fim

From the above lemma we obtain the following result.

\begin{proposicao}\label{From D to C} Let $\G$ be an ultragraph with no sinks that satisfy Condition~(RFUM) and let $M$ be a function from $\GG^1$ to $[0,1]$. Then there is a convex map between the set of functions $m:\GG^0\to [0,1]$ satisfying \ref{item:n1linha} to \ref{item:n4linha} of Proposition~\ref{From B to D} and the set of regular, Borel, probability measures $\mu$ on $X$ satisfying that $\mu(\theta_e(V))=M(e)\mu(V)$,
for all $e\in \GG^1$ and all Borel measurable subsets $V$ of $X_{e}^{-1}$.
\end{proposicao}
\demo

Using Carathéodory's extension theorem, there is a measure $\mu$ defined on the $\sigma$-algebra generated by $S$. Since $S$ forms a countable  basis, this is actually the Borel $\sigma$-algebra.

To prove that $\mu$ is probability, we first show that $\mu(D_{(A,A)})=m(A)$ for all $A\in \GG^0$. Using condition (RFUM), we can write $A=A_1\cup\cdots\cup A_n$, where $A_1,\ldots,A_{n-1}$ are minimal infinite emitters and $A_n$ is such that $\varepsilon(A_n)<\infty$. In particular, $D_{(A_i,A_i)}\in S$ for all $i=1,\ldots,n$. Using the inclusion-exclusion principle, we can write $\mu(A)$ as sum of $\pm\mu$ applied on $D_{(A_1,A_1)},\ldots,D_{(A_n,A_n)}$ and their intersections which are of the form $D_{(B,B)}$ and still belongs to $S$. By definition, $\mu(D_{(B,B)})=\kappa(D_{(B,B)})=m(B)$, whenever $D_{(B,B)}\in S$. Now, using the inclusion-exclusion principle again together with (m4') we conclude that $\mu(D_{(A,A)})=m(A)$.

Now, since $X=\bigcup_{A\in\GG^0}D_{(A,A)}$, using condition (m1') we see that $\mu(X)=1$. This also implies that $\mu$ is regular because all finite Borel measures are regular.

From equation (\ref{defkappa}), a convex combination of functions $m$ is preserved when passing to measures $\kappa$ on $S$, and therefore to measures $\mu$ on $X$.

That $\mu$ satisfies the scaling condition follows from the fact $\theta_e$ is bijective, so that it preserves unions and intersections, and that the condition holds for elements in the semi-ring $S$ as in Lemma \ref{measureonsemiring}.

\fim

\begin{proposicao}\label{equivalenciaEB} Let $\G$ be an ultragraph with no sinks that satisfy Condition~(RFUM). There is an affine bijection between $B^\beta$ and $E^\beta$. 
\end{proposicao}
\demo

Let $\psi$ be a state on ${C^*(\G)^\gamma}$ satisfying $
\psi(s_c p_A s_d^*)=\delta_{c,d}N(c)^{-\beta}\psi(p_{A\cap r(c)})$, for all $c \in \G^*$ and $A \in \GG^0$. Define 
a state $\omega$ on $C_0(X)$ by $\psi \circ \Phi^{-1}$, where $\Phi$ is the isomorphism of Theorem~\ref{crossedproduct} (see also Proposition~\ref{isoc0}).

We need to check the scaling condition. It is enough to check it in the dense subset of $C_0(X)$ given in Lemma~\ref{densealgebra} (see also Remark~\ref{rmkdensealgebra}). 

Let $e\in \G^1$. Notice that if $1_A \in C_0(X_{e^{-1}})$ then $A\subseteq X_{e^{-1}}=X_{r(e)}$. Then $$\omega(1_A \circ \theta_{e^{-1}}) = w(1_{A\cap r(e)} \circ \theta_{e^{-1}}) = \psi (s_e p_A s_e^*) = N(e)^{-\beta} \psi (p_{A\cap r(e)})$$$$\ \ \ \ \ \ \ \ \ = N(e)^{-\beta} \psi (\Phi^{-1}( 1_{A\cap r(e)})) = N(e)^{-\beta} \omega ( 1_{A\cap r(e)}).$$

Now let $c\in \G^*$ with $|c|\geq 1$, $A \in \G^0$, and suppose that $1_{r(c)\cap A} \circ \theta_{c^{-1}} \in C_0(X_{e^{-1}})$. Since $f\in C_0(X_{e^{-1}})$ we have that $s(c)\in r(e)$ and hence $r(c)= r(ec)$. Then
$$\omega(1_{r(c)\cap A} \circ \theta_{c^{-1}}\circ\theta_{e^{-1}}) = \omega(1_{r(ec)\cap A} \circ \theta_{(ec)^{-1}}) = \psi(s_{ec} p_A s_{ec}^*)= N(ec)^{-\beta}\psi (p_{A\cap r(ec)})$$
$$\ \ \ \ \ \ \ \ \ = N(e)^{-\beta}N(c)^{-\beta}\psi (p_{A\cap r(c)}) =N(e)^{-\beta} \psi (s_c p_A s_c^*)=N(e)^{-\beta} \omega(1_{r(c)\cap A} \circ \theta_{c^{-1}}).  $$
Hence the scaling condition is satisfied. 

To finish the proof, notice that the inverse of the correspondence above is given by the following map: Given a state $\omega$ on $C_0(X)$ that satisfies the scaling condition $\omega(f\circ\theta_e^{-1})=N(e)^{-\beta}\omega(f)$ for all $e\in \GG^1$ and all $f\in C_0(X_{e^{-1}})$, define $\psi$ as the state on ${C^*(\G)^\gamma}$ such that $\psi (s_c p_A s_d^*) = \delta_{c,d} \omega \circ \Phi$. Notice that $\psi\in E^\beta$, since for $c=c_1\ldots c_n \in \G^*$ and $A\in \G^0$, we have that
$$\begin{array}{lll}
\psi(s_c p_A s_c^*) & = & \omega(1_{r(c)\cap A} \circ \theta_{c^{-1}}) = N(c_1)^{-\beta} \omega(1_{r(c)\cap A} \circ \theta_{c_n^{-1}\ldots c_2^{-1}}) \\
 & = & N(c_1)^{-\beta}  N(c_2)^{-\beta} \omega(1_{r(c)\cap A} \circ \theta_{c_n^{-1}\ldots c_3^{-1}})= \ldots \\
 & = &  N(c_1)^{-\beta}  N(c_2)^{-\beta} \ldots N(c_n)^{-\beta} \omega(1_{r(c)\cap A}) = \\
  & = &  N(c)^{-\beta} \psi(p_{A\cap r(c)}).
 \end{array}$$

\fim

The next proposition is basically the same as \cite[Lemma~3.7~and~Proposition~3.8]{MR3539347}. We give the idea of the proof.

\begin{proposicao}\label{equivalenciaCB} Let $\G$ be an ultragraph with no sinks that satisfy Condition~(RFUM). There is an affine bijection between $B^\beta$ and $C^\beta$. 
\end{proposicao}
\demo

The Riesz representation Theorem gives the bijection between states $\omega$ on $C_0(X)$ and Borel regular probabilities $\mu$ on $X$. Restrict the corresponding measure $\mu$ to $X_{e^{-1}}$ and observe that $C_0(X_{e^{-1}})$ is dense in $L^1(X_{e^{-1}},\mu)$. If in one hand $\omega$ is in $B^\beta$, then the scaling condition holds for functions in $C_0(X_{e^{-1}})$ and then extends to elements of $L^1(X_{e^{-1}},\mu)$, in particular characteristic functions of measurable sets of $X_{e^{-1}}$. If on the other hand $\mu$ is in $C^\beta$, then the scaling condition holds for characteristic functions, which extends to elements of $L^1(X_{e^{-1}},\mu)$, and in particular to functions in $C_0(X_{e^{-1}})$.

\fim

\begin{remark} We end the section with the proof of Theorem~\ref{teoremaprincipal}.
\end{remark}
\demo

Notice that from Proposition~\ref{equivalenciaAE} we have an affine isomorphism between $A^\beta$ and $E^\beta$), from Proposition~\ref{equivalenciaEB} we have an affine isomorphism between $E^\beta$ and $B^\beta$), and from Proposition~\ref{equivalenciaCB} we have an affine isomorphism between $C^\beta$ and $B^\beta$). 

Taking $M:\G^1\rightarrow [0,1]$ defined by $M(e)=(N(e))^{-\beta}$ on Proposition~\ref{From B to D} we obtain an affine map from $B^\beta$ to $D^\beta$, and taking the same $M$ on Proposition~\ref{From D to C} we obtain an affine map from $D^\beta$ to $C^\beta$. The result now follows from the fact that the maps from $B^\beta$ to $D^\beta$, from $D^\beta$ to $C^\beta$ and from $C^\beta$ to $B^\beta$ compose to the identity.

\fim

\section{Ground States}\label{ground}
In this section we apply some of our previous results to characterize the set of ground states on the C*-algebra of an ultragraph with no sinks that satisfy Condition~(RFUM) (we recall the general definition of ground states in the beginning of Section~\ref{KMSgeneral}).

Let $\G$ be an ultragraph with no sinks that satisfy Condition~(RFUM). Define 

\begin{itemize}
		\item[$A^{gr}$:] the set of ground states in $C^*(\G)$,
		\item[$B^{gr}$:] the set of states $\omega$ of $C_0(X)$ such that $\omega(1_e)=0 $, for all $e\in \GG^1$,
	   \item[$C^{gr}$:] the set of regular Borel probability measures $\mu$ on $X$  such that $\mu(A) = 0 $, for every $e\in \GG^1$ and every Borel measurable subset $A$ of $X_{e}$,
	  \item[$D^{gr}$:] the set of functions $m:\GG^0\to [0,1]$ satisfying 
\begin{enumerate}\renewcommand{\theenumi}{m\arabic{enumi}}
  	\item\label{item:gn1} $\lim_{A\in \GG^0} m(A) =1$;
    \item\label{item:gn2} $m(A)=0 $ if $|\varepsilon(A)|<\infty$;
    \item\label{item:gn3} $m(A\cup B) = m(A)+m(B)-m(A\cap B)$.
  \end{enumerate}
	\end{itemize}

Next we prove the following:

\begin{teorema}\label{teoremaground} Let $\G$ be an ultragraph with no sinks that satisfy Condition~(RFUM). Then there exists an affine isomorphism between $A^{gr}$, $B^{gr}$, $C^{gr}$, and $D^{gr}$.

\end{teorema}
\demo

The existence of an affine isomorphism between $A^{gr}$ and the set of states $\phi$ of $C_0(X)$ such that $\phi(f)=0 $, for all $e\in \GG^1$ and $f\in C(X_e)$ follows from Theorem~4.3 of \cite{MR1953065}. Since $\phi$ is a state, and $1_e$ is a unit for $C(X_e)$, it follows that if $\phi(1_e)=0$ then $\phi(f)=0 $ for all $f\in C(X_e)$. Therefore we have that $A^{gr}$ is isomorphic to $B^{gr}$, via an affine isomorphism.

As with KMS states, an affine isomorphism between $B^{gr}$ and $C^{gr}$ is obtained analogously to Proposition~4.8 in \cite{MR3539347}.

Finally, an affine isomorphism between $B^{gr}$ to $D^{gr}$ is obtained by application of Propositions~\ref{From B to D} and \ref{From D to C} with $M(e)=0$ for all $e\in \G^1$, and proceeding as in the proof of Theorem~\ref{teoremaprincipal}.
\fim

\section{Example}\label{secex}

Below we recall the construction, given in \cite{MR2001938}, of an ultragraph such that the associated C*-algebra is neither a graph C*-algebra nor an Exel-Laca algebra.

\begin{definicao}
	If $I$ is a countable set and $A$ is an $I
	\times I$ matrix with entries in $\{ 0, 1\}$, then we may
	form the ultragraph $\G_A := (G_A^0, \G_A^1, r, s)$
	defined by $G_A^0 := \{v_i : i \in I \}$, $\G_A^1 :=  \{e_i : i \in I \}$,
	$s(e_i) = v_i$ for all $i \in I$, and $r(e_i)=\{v_j : A(i,j)
	= 1 \}$.  
	\label{edgeLG}
\end{definicao}

Let $A$ be the countably infinite matrix
$$
A = \left(\begin{smallmatrix}
1 & 0 & 0 & 1 & 1 & 1 & 1 &\\
0 & 1 & 0 & 1 & 1 & 1 & 1 &\\
0 & 0 & 1 & 1 & 1 & 1 & 1 &\\
1 & 0 & 0 & 1 & 0 & 0 & 0 & \cdots \\ 
0 & 1 & 0 & 0 & 1 & 0 & 0 & \\ 
0 & 0 & 1 & 0 & 0 & 1 & 0 & \\ 
0 & 0 & 0 & 1 & 0 & 0 & 1 & \\
&   &   & \vdots & & & & \ddots \\ 
\end{smallmatrix}\right).
$$

For the matrix $A$ above, let $\G := (G^0, \G^1, r, s
)$ be the ultragraph $\G_A$ of Definition~\ref{edgeLG}. 
We define an ultragraph $\FF$ by adding a single vertex $\{
w \}$ to $\G$ and a countable number of edges with source $w$
and range $G^0$.  More precisely, we define $\FF := (F^0,
\FF^1, r, s)$ by
$$ F^0 := \{w \} \cup G^0 \quad \quad \quad \quad \FF^1 :=
\{f_i\}_{i=1}^\infty \cup \G^1 $$ and we extend $r$ and
$s$ to $\FF^1$ by defining $s(f_i) = \{w \}$ and $r(f_i) =
G^0$ for all $1 \leq i < \infty$.

Our goal is to define a suitable function $N$ for which we can find $m$ that satisfies conditions m1 to m4 defining the set $D^\beta$. Using Definition~\ref{edgeLG} for the matrix $A$ above, we conclude that the set $\FF^0$ is composed by all finite sets of vertices, and finite sets union with the set $B:=\{v_4,v_5,\ldots\}$. In particular $F^0$ is an element of $\FF^0$.

Now, if a function $m$ satisfying m1 to m4 is to exist, for a given $N$ and $\beta$, using m1 and m4, we must have that
\begin{equation}\label{eq.example}
	m(\{w\})+m(\{v_1\})+m(\{v_2\})+m(\{v_3\})+m(B)=m(F^0)=1.
\end{equation}
Also, for $i=1,2,3$, using conditions m2 and m4, it should hold that 
\[m(\{v_i\})=N(e_i)^{-\beta}(m(\{v_i\})+m(B)),\]
so that we can describe $m(\{v_i\})$, for $i=1,2,3$, depending on $N$, $\beta$ and $m(B)$ as
\[m(\{v_i\})=\frac{N(e_i)^{-\beta}}{1-N(e_i)^{-\beta}}m(B).\]

For $i\geq 4$, we have that $r(e_i)=\{v_{i-3},v_i\}$. Using an induction on each congruence modulo 3, we see that $m(\{v_i\})$ can be written depending on $N$, $\beta$ and $m(B)$ using the edges $e_j$ for $j$ in the congruence class modulo 3 up to $i$. For simplicity, let us assume that for all $1 \leq i < \infty$, $N(e_i)=d$ for a fixed number $d\in(1,\infty)$ and define $d_{\beta}=\frac{d^{-\beta}}{1-d^{-\beta}}$. Notice that $d_{\beta}$ goes to $0$ as $\beta$ goes to infinity. We claim that for $1 \leq i < \infty$, writing $i=3q+r$ with $r=1,2,3$, we have that
\begin{equation}\label{eq.vi}
	m(\{v_i\})=d_{\beta}^{q+1}m(B).
\end{equation}
This is the case for $q=0$. Now suppose it is true for a fixed $r=1,2,3$ and a given $q\geq 0$. Then
\[m(v_{3(q+1)+r})=d^{-\beta}(m(v_{3(q+1)+r})+m(v_{3q+r}))\]
and hence
\[m(v_{3(q+1)+r})=d_{\beta}m(v_{3q+r})=d_{\beta}^{q+2}m(B).\]

Using this one can show that, if
\begin{equation}\label{eq.exemaple.B}
	\frac{6d_\beta^2}{1-d_\beta^2}\leq 1,
\end{equation}
then m3 holds for $B$ and any finite subset $F$ of $\varepsilon(B)$.

Now, choose a sequence $\{c_i\}_{i=1}^{\infty}$ with $c_i\in(1,\infty)$ and define $N(f_i)=c_i$. Using equation (\ref{eq.example}), we can show that for m3 to be true for $\{w\}$, it is necessary that for all $k\in\mathbb{N}$, 
\[\frac{m(\{w\})}{1-m(\{w\})}\geq \sum_{i=1}^k c_i^{-\beta}.\]

It follows that if $\beta$ is such that the series $\sum_{i=1}^{\infty} c_i^{-\beta}$ converges and that (\ref{eq.exemaple.B}) holds, then for any value of $m(\{w\})\in[0,1)$ such that
\[\frac{m(\{w\})}{1-m(\{w\})}\geq \sum_{i=1}^{\infty} c_i^{-\beta},\]
we can define $m(B)$ and $m(\{v_i\})$ for $i\in \mathbb{N}$ using equations (\ref{eq.vi}) and (\ref{eq.exemaple.B}).

Also, independently of $N$ and $\beta$, the choice $m(\{w\})=1$ and $m(B)=m(\{v_i\})=0$ for $i\in \mathbb{N}$ always satisfies m1-m4.

In the case of ground states, a function $m:\GG^0\rightarrow [0,1]$ is is $D^{gr}$ if, and only if,  $m(\{v_i\})=0$ for $i\in \mathbb{N}$, and $m(B)+m(\{w\})=1$.

\section{Acknowledgements}

The authors would like to thank Zahra Afsar for valuable discussions regarding the present paper. In particular, the second author would like to thank Zahra for teaching him the theory of KMS states.

\bibliographystyle{abbrv}
\bibliography{references}

\vspace{1.5pc}

Gilles Gonçalves de Castro, Departamento de Matemática, Universidade Federal de Santa Catarina, Florianópolis, 88040-900, Brazil.

Email: gilles.castro@ufsc.br	

\vspace{0.5pc}

Daniel Gonçalves, Departamento de Matemática, Universidade Federal de Santa Catarina, Florianópolis, 88040-900, Brazil.

Email: daemig@gmail.com

\end{document}